\documentclass[11pt]{amsart}
\usepackage[latin1]{inputenc}
\usepackage{amsmath,amsthm, amscd, amssymb, amsfonts}
\usepackage[dvips]{graphicx}
\numberwithin{equation}{section}

\theoremstyle{plain}

\newtheorem{teor}{Theorem}[section]
\newtheorem{corollary}[teor]{Corollary}
\newtheorem{proposition}[teor]{Proposition}
\newtheorem{lemma}[teor]{Lemma}

\theoremstyle{definition}

\newtheorem{examples}[teor]{Examples}

\theoremstyle{remark}
\newtheorem{remark}[teor]{Remark}

\newcommand\tright{\vartriangleright}
\newcommand\tleft{\vartriangleleft}
\newcommand\id{\operatorname{id}}
\newcommand\ad{\operatorname{ad}}

\newcommand\Bigal{\operatorname{Bigal}}

\newcommand\Hom{\operatorname{Hom}}
\newcommand\Rep{\operatorname{Rep}}
\newcommand\Ind{\operatorname{Ind}}
\newcommand\Opext{\operatorname{Opext}}

\newcommand\vect{\operatorname{Vec}}
\newcommand\Aut{\operatorname{Aut }}
\newcommand\C{\mathcal C}

\begin{document}

\title[Hopf algebra extensions of group algebras]{Hopf algebra extensions of group algebras and Tambara-Yamagami categories}
\author{Sonia Natale}
\address{Facultad de Matem\'atica, Astronom\'\i a y F\'\i sica,
Universidad Nacional de C\'ordoba, CIEM -- CONICET, (5000) Ciudad
Universitaria, C\'ordoba, Argentina} \email{natale@mate.uncor.edu}
\thanks{This work was partially supported by CONICET, Fundaci\' on Antorchas, ANPCyT-Foncyt and Secyt (UNC)}
\date{Revised version of April 2009}

\begin{abstract} We determine the structure of Hopf algebras that admit an
extension of a group algebra by the cyclic group of order 2. We
study the corepresentation theory of such Hopf algebras, which
provide a generalization, at the Hopf algebra level, of the so
called Tambara-Yamagami fusion categories.  As a byproduct, we
show that every semisimple Hopf algebra of dimension $<36$ is
necessarily group-theoretical; thus $36$ is the smallest possible
dimension where a non group-theoretical example occurs.
\end{abstract}

\maketitle
\section{Introduction}

Along this paper we shall work over an algebraically closed base
field $k$ of characteristic zero.

An important problem related to the classification of semisimple
Hopf algebras over $k$ was the question raised in the paper
\cite{ENO}, whether all semisimple Hopf algebras over $k$ are
group-theoretical; see \cite[Section 8]{ENO} for a discussion of
group-theoretical categories and quasi-Hopf algebras. A positive
answer to this question was obtained for certain cases in the
paper \cite{dgno}; in particular, semisimple Hopf algebras of
dimension $p^n$, where $p$ is a prime number, are always
group-theoretical.

In the general case, the question was answered negatively by
Nikshych \cite{nik}. In fact, a family of semisimple Hopf algebras
which are not group-theoretical is constructed in \cite{nik}: a
Hopf algebra $H$ in this family fits into an exact sequence
\begin{equation}\label{ext-nik}k \to k^{\mathbb Z_2} \to H \to
(kG)^J \to k,\end{equation} where $G$ is a certain finite group
and $J \in kG\otimes kG$ is an invertible twist. It turns out that
these examples are all \emph{semisolvable} in the sense of
\cite{MW}.

\medbreak In this paper we consider semisimple Hopf algebras which
fit into an exact sequence \begin{equation}\label{ext-triang}k \to
A \to H \to k{\mathbb Z_2} \to k,\end{equation} where $A$ is a
triangular semisimple Hopf algebra. In view of the classification
of semisimple triangular Hopf algebras, $A \simeq (kG)^J$ as Hopf
algebras, for some finite group $G$ and invertible twist $J \in
kG\otimes kG$ \cite{eg-triang}.

Contrary to the situation for the extensions \eqref{ext-nik},
those in \eqref{ext-triang} are always group-theoretical, despite
of the symmetry in the form of the extensions. Moreover, we show
in Theorem \ref{triangular} that a Hopf algebra $H$ as in
\eqref{ext-triang} is twist equivalent to an \emph{abelian}
extension, and \textit{a fortiori} group-theoretical by the main
result of \cite{gp-ttic}.

As mentioned before, a Hopf algebra $H$ as in \eqref{ext-triang}
is twist equivalent to a semisimple Hopf algebra whose group of
group-like elements has index 2. Therefore, up to twist
deformations, we are lead to considering extensions of the form
\begin{equation}\label{ext-grupo}k \to kG \to H \to k{\mathbb Z_2}
\to k,\end{equation} where $G$ is a finite group. We prove that if
such a nontrivial $H$ exists, then $H$ fits also into an abelian
(cocentral) exact sequence; in particular, the group $G \simeq
G(H)$ is necessarily an extension of an abelian group of order
$d^2$, $d
> 1$; see Propositions \ref{abeliana}, \ref{f0}.

In the situation of \eqref{ext-grupo}, we show that there is an
equivalence of fusion categories $\Rep H \simeq (\Rep G)^{\mathbb
Z_2}$, where $(\Rep G)^{\mathbb Z_2}$ is an appropriate $\mathbb
Z_2$-equivariantization of $\Rep G$. More generally, we show that
the representation category of every cocentral Hopf algebra
extension is an equivariantization; see Proposition
\ref{equiv-hopf}.

On the other hand,  the corepresentation category of a Hopf
algebra $H$ as in \eqref{ext-grupo} turns out to be a
generalization, at the Hopf algebra level, of the so called
Tambara-Yamagami  categories \cite{TY}: these are fusion
categories with isomorphism classes of simple objects
parameterized by the set $\Gamma \cup \{ x\}$, where $\Gamma$ is a
finite (necessarily abelian) group, $x \notin \Gamma$, obeying the
fusion rules
$$s\otimes t = st, \quad s, t \in \Gamma, \qquad x \otimes x = \oplus_{s \in \Gamma}s.$$
Distinct features of Tambara-Yamagami fusion rules and some of
their generalizations have attracted the interest of several
authors in the last years, c.f., \cite{artamonov, ENO, nik,
siehler1, siehler2}.

In the situation of \eqref{ext-grupo} there may be more than one
non-invertible object, nevertheless, all non-invertible objects
have the same dimension $d$ and the product of any two of them is
a sum of invertible objects belonging to a fixed normal abelian
subgroup of order $d^2$ of $G$.

\medbreak We prove the following classification theorem.

\begin{teor}\label{cls} Semisimple Hopf algebras $H$ with $[H: kG(H)] =
2$ are determined by triples $(\Gamma, F, \xi)$, where

\begin{enumerate}\item\label{1} $\Gamma$ is a finite abelian group,
\item\label{3} $F$ is a finite group acting on $\Gamma$ by group automorphisms,
\item\label{2} $\xi$ is an element of the
group $\Opext(k^{\widehat \Gamma}, kF)$, \end{enumerate}
satisfying
\begin{flalign}&[F: F_0] = 2, \text{where } F_0 \text{ is the subgroup } F_0 : = \{ x \in F:\, [\tau_x] = 1
\}, \text{and} & \\ & \tau_x \text{ is a non-degenerate 2-cocycle,
for all } x\in F\backslash F_0,& \end{flalign} if $(\sigma, \tau)$
is a pair of compatible cocycles representing $\xi$, where $\tau_x
\in Z^2(\widehat \Gamma, k^{\times})$ is defined by $\tau_x(s, t)
= \tau(s, t) (x)$, $x \in F$, $s, t \in \widehat \Gamma$.

If $H$ corresponds to the triple $(\Gamma, F, \xi)$, then $G(H)$
is isomorphic to the crossed product $\Gamma \rtimes_{\sigma}
F_0$. \end{teor}

Theorem \ref{cls} is proved in Section \ref{clsfusionrules}. It
extends the classification result \cite[Theorem 3.5]{T} of
Tambara; see Remark \ref{cls-tambara}.

\medbreak  The smallest example of a non group-theoretical
semisimple Hopf algebra from the construction in \cite{nik} has
dimension $36$ and it is a semisolvable Hopf algebra.

We also show in this paper that every semisimple Hopf algebra of
dimension $<36$ is group-theoretical. So that, in fact, $36$ is
the smallest possible dimension that a non group-theoretical
semisimple Hopf algebra can have. See Theorem \ref{menor36}. It is
known that in dimension $< 36$ every semisimple Hopf algebra is
upper and lower semisolvable; see \cite{Na1}. Coincidentally, the
first non semisolvable example also appears  in dimension $36$
\cite{twist-simple}.

\medbreak Moreover, in dimension $<36$, except for dimension $24$,
every semisimple Hopf algebra is either nilpotent (dimensions $p$,
$p^2$, $p^3$, $p^4$, $p^5$, where $p$ is a prime number) or an
abelian extension (dimensions $30$, $pq$, and $pq^2$, where $p$
and $q$ are distinct prime numbers). Then all these Hopf algebras
are group-theoretical, in view of \cite{dgno, gp-ttic}. Therefore
we only need to consider the case of dimension $24$. We prove
that, in fact, up to a cocycle twist of the multiplication or the
comultiplication, every semisimple Hopf algebra of dimension $24$
fits into an abelian exact sequence; see Proposition \ref{main}.

\medbreak The paper is organized as follows. In Section
\ref{prels} we recall some preliminary notions on Hopf algebra
extensions needed later on. We also discuss in Subsection
\ref{gpttic} some properties of group-theoretical Hopf algebras
and their relation with Hopf algebra extensions. In Section
\ref{repsbygroup} we describe the representation theory of
cocentral extensions of a Hopf algebra. In Sections
\ref{fusionrules} and \ref{clsfusionrules} we prove our main
results on Hopf algebras with group of group-likes of index $2$.
In Section \ref{lowdim} we apply the above to prove our statement
on semisimple Hopf algebras of dimension $<36$. At the end of the
paper we include an appendix where we give a sufficient condition
for normality of a group-like Hopf subalgebra in a semisimple Hopf
algebra.

\medbreak \noindent \textbf{Conventions and notation.} The
notation for Hopf algebras is standard: $\Delta$, $\epsilon$, $S$,
denote the comultiplication, counit and antipode, respectively. We
refer the reader to \cite{ENO, kassel} for the terminology and
notation on tensor categories and tensor functors used throughout.

Let $d$ be a positive integer. The notation $M_d(k)$ will indicate
the matrix coalgebra of dimension $d^2$. For a group $G$, the
group algebra of $G$ and its dual will be denoted by $kG$ and
$k^G$, respectively. The group of one-dimensional characters of
$G$ will be indicated by $\widehat G$.

\medbreak \noindent \textbf{Acknowledgement.} The author thanks
Gast\' on Garc\' \i a for discussions on the contents of Section
\ref{lowdim}, that motivated the inclusion of that section in this
paper.

\section{Hopf algebra extensions}\label{prels}

Let $H$ be a Hopf algebra over $k$. The left (respectively, right)
adjoint action of $H$ on itself is defined by $\ad_h (x)= h_1 x
\mathcal S (h_2)$ (respectively, $\ad^r_h (x)= \mathcal S (h_1) x
h_2$), $x, h\in H$. A Hopf subalgebra $K \subseteq H$ is called
\emph{normal} if it is stable under both adjoint actions; $H$ is
called \emph{simple} if it contains no proper normal Hopf
subalgebras.

Let $K \subseteq H$ be a normal Hopf subalgebra. Then $B = H/HK^+$
is a quotient Hopf algebra and the sequence of Hopf algebra maps
$k \to K \to H \to B \to k$ is an exact sequence of Hopf algebras.
In this case $H$ is called an \emph{extension} of $K$ by $B$. If
the extension is cleft, then  $H$ is isomorphic to a bicrossed
product $H \cong K {}^{\tau}\#_{\sigma} B$ as a Hopf algebra. See
\cite{AD, ma-ext, ma-ext2} and references therein.

Suppose $H$ is finite dimensional. Then  $K$ is normal in $H$ if
and only if $K$ is stable under the left adjoint action. In this
case, the extension $K\subseteq H$ is always cleft. Moreover, in
the finite dimensional context, the notion of simplicity is
self-dual; that is, $H$ is simple if and only if $H^*$ is.

\medbreak By \cite[Corollary 1.4.3]{Na1}, if $K \subseteq H$ is a
normal Hopf subalgebra, such that $\dim K$ is the least prime
number dividing $\dim H$, then $K$ is central in $H$.

\subsection{Abelian extensions}\label{ab-gral}
We refer the reader to \cite{ma-ext, ma-ext2} for the notion of
abelian exact sequence and, in particular, for the study of the
cohomology theory underlying such exact sequence.

Suppose that $L = F\Gamma$ is an exact factorization of the finite
group $L$, where $\Gamma$ and $F$ are subgroups of $L$.
Equivalently, $F$ and $\Gamma$ form a matched pair of finite
groups with the actions $\vartriangleleft : \Gamma \times F \to
\Gamma$, $\vartriangleright : \Gamma \times F \to F$, defined by
$sx = (x \vartriangleleft s)(x \vartriangleright s)$, $x \in F$,
$s \in \Gamma$.

Let $\sigma: F \times F \to (k^\Gamma)^{\times}$, $\sigma(x, y) =
\sum_s\sigma_s(x, y)e_s$, and $\tau: \Gamma \times \Gamma \to
(k^F)^{\times}$, $\tau(s, t) = \sum_x\tau_x(s, t)e_x$, be
normalized 2-cocycles with the respect to the actions afforded,
respectively, by $\vartriangleleft$ and $\vartriangleright$,
subject to appropriate compatibility conditions \cite{ma-ext}.
Here, $e_y \in k^F$,  $y \in F$,  are the canonical idempotents
defined by $e_y(x) = \delta_{x, y}$, and similarly for $e_s \in
k^{\Gamma}$.

Consider the bicrossed product $H = k^\Gamma \,
{}^{\tau}\#_{\sigma}kF$ corresponding to this data. Then $H$ is a
Hopf algebra, with multiplication and comultiplication determined
by
\begin{align}\label{mult}(e_s \# x)  (e_t \# y)
& =  \delta_{s \vartriangleleft x, t} \, \sigma_s(x, y)  \, e_s \# xy, \\
\label{comult}\Delta(e_s \# x) & = \sum_{gh = s} \tau_x(g, h) \,
e_g \# (h \vartriangleright x) \otimes e_h \# x,
\end{align}
$s, t \in \Gamma$, $x, y \in F$,
 and there is an \emph{abelian} exact sequence $k \to
k^{\Gamma} \to H \to kF \to k$. Moreover, every Hopf algebra
fitting into such exact sequence can be described in this way.
This gives a bijective correspondence between the equivalence
classes of these Hopf algebra extensions and a certain abelian
group $\Opext (k^\Gamma, kF)$ associated to the matched pair $(F,
\Gamma)$.

\begin{remark}\label{cocentralabeliana}  Suppose that the extension
$k \to k^{\Gamma} \to H \to kF \to k$ is cocentral or,
equivalently, that the action $\vartriangleright : \Gamma \times F
\to F$ is trivial, so that $L = F\Gamma \simeq \Gamma \rtimes F$
is a semidirect product. By Proposition \ref{equiv-hopf} below,
there is an equivalence of tensor categories $$\Rep H \simeq
(\vect^{\Gamma})^F,$$ where $\vect^{\Gamma}$ denotes the tensor
category of $\Gamma$-graded vector spaces and the $F$-action is
given by $T_x(V)_s = V_{s \vartriangleleft x}$, $s \in \Gamma$,
$x\in F$.
\end{remark}

\medbreak The following lemma describes the group of group-likes
in an abelian extension. Let $F^\Gamma \subseteq F$ denote the
(subgroup) of elements in $F$ which are invariant under the action
$\vartriangleright$.

Note that if $x \in F^\Gamma$, then $\tau_x: \Gamma \times \Gamma
\to k^{\times}$ is a normalized 2-cocycle. We shall denote by
$[\tau_x]$ its cohomology class in $H^2(\Gamma, k^{\times})$.

\begin{lemma}\label{grouplikes} The exact sequence of Hopf algebras $k \to
k^{\Gamma} \to H \to kF \to k$ induces by restriction an exact
sequence of groups $$1 \to \widehat \Gamma \to G(H) \to F_0 \to
1,$$ where $F_0 = \{ x \in F^\Gamma: \, [\tau_x] = 1 \}$.
\end{lemma}

Therefore $G(H)$ is isomorphic to a crossed product $G(H) \simeq
\widehat \Gamma \rtimes_{\vartriangleleft, \sigma} F_0$.

\begin{proof} The proof is straightforward using the formula for the comultiplication \eqref{comult}.
Indeed, each group-like element of $H$ is of the form $\gamma \#
x$, where $0 \neq \gamma \in k^{\Gamma}$ and $x \in F^{\Gamma}$,
such that $\gamma (1) = 1$ and $\tau_x(s, t) =
\gamma(s)\gamma(t)\gamma(st)^{-1}$, $\forall s, t \in \Gamma$.
\end{proof}

\begin{remark}\label{estab} Formula \eqref{comult} implies that for all $x\in F$,
the subspace $I_x = k^{\Gamma}\# x \subseteq H$ is a left coideal,
and we have a decomposition $H = \oplus_{x \in F}I_x$. Moreover,
by \eqref{mult}, $I_xI_y \subseteq I_{xy}$, $\forall x, y \in F$.

Suppose $x \in F^{\Gamma}$. Then $I_x$ is a subcoalgebra of $H$,
and $I_x = (k_{\tau_x}\Gamma)^*\# x$, as coalgebras. In
particular, $I_x \simeq k^{\Gamma}$ if and only if $[\tau_x] = 1$
and it is simple if and only if $\tau_x$ is non-degenerate.

If $y \in F^{\Gamma}$ is such that $\tau_y$ is non-degenerate,
then the relation $I_1I_y \subseteq I_{y}$ implies that $\widehat
\Gamma = G(k^{\Gamma}) \subseteq G(H)$ stabilizes the simple
subcoalgebra $I_y$. \end{remark}

\subsection{Relations with group-theoretical Hopf
algebras}\label{gpttic}

Let $G$ be a finite group and let $\omega$ be a 3-cocycle on $G$.
Let $\C = \C(G, \omega)$ be the category of $G$-graded vector
spaces with associativity isomorphism given by $\omega$. Let also
$F$ be a subgroup of $G$ and $\alpha$ a 2-cochain on $F$. Suppose
that  $\omega\vert_F = d\alpha$; so that the twisted group algebra
$k_{\alpha}F$ is an algebra in $\C$. Then the category $\C(G,
\omega, F, \alpha)$ of $k_{\alpha}F$-bimodules in $\C$ is a fusion
category, called a \emph{group-theoretical} category. A
(quasi-)Hopf algebra is called group-theoretical if its category
of representations is group-theoretical \cite[8.8]{ENO}.

\medbreak In particular, the class of group-theoretical
(quasi-)Hopf algebras is closed under twisting. Furthermore, it is
shown in \cite{ENO} that the class of group-theoretical categories
is closed under duals, Drinfeld centers, taking full fusion
subcategories and components in a quotient category.

\medbreak It was shown in \cite{gp-ttic} that every abelian
extension of Hopf algebras is group-theoretical. However, the
group-theoretical Hopf algebras are not closed under taking
extensions. In the paper \cite{nik}, Nikshych constructs examples
of semisimple Hopf algebras $H$ which are not group-theoretical,
and nevertheless fit into an exact sequence $k \to k^{\mathbb Z_2}
\to H \to (kG)^J \to k$, where $G$ is a certain finite group and
$J \in kG\otimes kG$ is an invertible twist.

\section{Representation category of a cocentral extension}\label{repsbygroup}

\subsection{$G$-equivariantization}

We start this subsection by recalling  the definition of
equivariantization of a fusion category. This notion has been
considered by different authors, see  \cite{agaitsgory, fw, nik,
tambara}. We first present the construction for group actions on a
linear category $\mathcal C$, that is, dropping the tensor
structure in $\mathcal C$.

Let $F$ be a finite group and let $\mathcal C$ be a $k$-linear
category. The group $F$ will also be regarded as a monoidal
category, denoted by $\underline F$, whose objects are the
elements of $F$, arrows are identities and tensor product is the
multiplication in $F$.

\medbreak Let $\underline \Aut\, \mathcal C$ denote the monoidal
category whose objects are autoequivalences of $\mathcal C$,
morphisms are isomorphisms of functors and tensor product is given
by composition of functors. An action of $F$ on $\mathcal C$ is a
monoidal functor
\begin{equation}\label{action}(T, f): \underline F \to \underline \Aut\, \mathcal
C.\end{equation}

In particular, for all $g, h \in F$, there are natural
equivalences
\begin{equation}f_{g, h}: T_g \circ T_h \to T_{gh},
\end{equation} giving the tensor structure to the functor $T$, where $T_g
= T(g)$, $g \in F$.

\medbreak Suppose that \eqref{action} is an action of $F$ on
$\mathcal C$. An $F$-equivariant object in $\mathcal C$ is a pair
$(V, (u_g^V)_{g \in F})$, where $V$ is an object of $\mathcal C$
and $u^V_g: T_g(V) \to V$, $g \in F$, are isomorphisms compatible
with the tensor structure on $T$ in the sense that, for all $g, h
\in F$,
\begin{equation}\label{deltau} u^V_g T_g(u^V_h) = u^V_{gh} f^V_{g, h}.\end{equation}
An $F$-equivariant morphism $\phi: (U, u_g^U) \to (V, u_g^V)$
between $F$-equivariant objects $(U, u_g^U)$ and $(V, u_g^V)$, is
a morphism $\phi: U \to V$ in $\mathcal C$ such that $\phi u^U_g =
u^V_g\phi$, for all $g \in F$.

The $F$-\emph{equivariantization} of $\mathcal C$, denoted
$\mathcal C^F$,  is defined to be the category of $F$-equivariant
objects and $F$-equivariant morphisms.

\medbreak If $\mathcal C$ is a tensor (respectively, fusion)
category, and $F$ acts by tensor autoequivalences, that is,  the
action \eqref{action} takes values in the subcategory $\underline
\Aut_{\otimes}\mathcal C$  of tensor autoequivalences of $\mathcal
C$ and isomorphisms of tensor functors, then the
equivariantization $\mathcal C^F$ is a tensor (respectively,
fusion) category with tensor product inherited from $\mathcal C$.

More precisely,  it is required in this case that $(T_g, j_g)$ be
a monoidal functor, with  $j_g\vert_{U, V}: T_g(U\otimes V)
\overset{\simeq}\to T_g(U)\otimes T_g(V)$,  for all $g \in F$.
Then, for equivariant objects $(U, u_g^U)$ and $(V, u^V_g)$, their
tensor product $(U\otimes V, u_g^{U\otimes V})$ is an equivariant
object, where $u_g^{U\otimes V} = (u_g^U \otimes
u_g^V)j_g\vert_{U, V}$. Moreover, $\mathcal C^F$ is dual to a
crossed product fusion category $\mathcal C \rtimes F$ with
respect to the indecomposable module category $\mathcal C$. See
\cite{nik}.

\begin{remark} Suppose $\mathcal C = {}^H\mathcal M$ is the
category of (finite dimensional) left $H$-comodules over the Hopf
algebra $H$. If $L$ is another Hopf algebra, it is known that
isomorphism classes of equivalences of tensor categories
${}^H\mathcal M \to {}^L\mathcal M$ are in bijective
correspondence with classes of $(L, H)$-bigalois extensions
\cite{scha-bigalois}.

Let $\Bigal (H)$ denote the set of classes of $(H, H)$-bigalois
extensions. Then $\Bigal (H)$ is a group under cotensor product,
and there is an inclusion of tensor categories $\underline{\Bigal
(H)} \to \underline \Aut_{\otimes}\mathcal C$. Moreover, every
object of $\Aut_{\otimes}\mathcal C$ is isomorphic to exactly one
object of $\underline{\Bigal (H)}$.
\end{remark}

\subsection{Representations of crossed
products}\label{repcrossprod}
 Let $A$ be a finite dimensional $k$-algebra. Suppose
that the finite group $F$ measures $A$ via $.:kF \otimes A \to A$,
$g \otimes a \mapsto g.a$ and let  $\sigma: kF \otimes kF \to A$
be an invertible compatible 2-cocycle; that is, $\sigma(g, h)$ is
invertible in $A$,  and we have
\begin{align}& g.(ab) = (g.a)(g.b), \quad g.1 = 1, \\
& (g. \sigma(h, t))\sigma(g, ht) = \sigma(gh, t)\sigma(g, h),
\quad \sigma(1, g) = \sigma(g, 1) = 1, \\
\label{asoc}& g.(h.a) = \sigma(g, h) (gh.a)\sigma(g, h)^{-1},
\end{align}
for all $g, h, t \in F$, $a, b \in A$. See \cite[Chapter 7]{Mo}.
In this situation, we shall also say that $.$ is a weak action.
Let $A \#_{\sigma}kF$ be the corresponding crossed product.

\medbreak The category $\Rep A$ of finite dimensional
representations of $A$ admits an action of the group $F$, $(T, f):
\underline F \to \underline \Aut(\Rep A)$,  defined as follows:
$T_g(V) = V$, with $a._gv = (g.a)v$, $V \in \Rep A$, $v \in V$. We
let $T_g = \id$ on morphisms,  and
$$f_{g, h}: T_g(T_h(V)) \to T_{gh}(V), \quad v \mapsto \sigma(g,
h)^{-1}v,$$ for all $g, h \in F$. Note that $f_{g, h}$ is an
$A$-linear isomorphism  thanks to \eqref{asoc}.

\medbreak Suppose $V$ is a $F$-equivariant object with respect to
this action. This means that there are $A$-linear isomorphisms
$u_g = u_g^V: V \to V$ satisfying
\begin{equation}\label{g-act1}(g.a)|_V =
u_g^{-1}a|_Vu_g,\end{equation} for all $g \in F$, $a \in A$. By
\eqref{deltau}, we have in addition
\begin{equation}\label{g-act2}\sigma(g, h)^{-1}|_V = f_{g, h}|_V =
u_{gh}^{-1}u_gu_h,\end{equation} for all $g, h \in F$.

It is straightforward to check that conditions \eqref{g-act1} and
\eqref{g-act2} imply that there is a well-defined action of the
crossed product $A\#_{\sigma}kF$ on $V$ determined by
\begin{equation} (a\#g).v = au_g^{-1}(v),
\end{equation} for $a \in A$, $g \in F$, $v \in V$. Moreover,
morphisms of $F$-equivariant objects are exactly morphisms of
$A$-modules commuting with the action of $F$ afforded by the
$u_g$'s, so they are $A\#_{\sigma}kF$-morphisms. We get in this
way a functor
\begin{equation} \mathcal F: (\Rep A)^F \to \Rep (A\#_{\sigma}F).
\end{equation}

\begin{proposition}\label{equiv} The functor $\mathcal F: (\Rep A)^F \to \Rep H$ defines an equivalence
of  categories. \end{proposition}

\begin{proof} Suppose that $W \in \Rep (A\#_{\sigma}F)$. Then $W$
is a representation of $A$ by restriction. Moreover,  $(W, u^W_g)$
is a $F$-equivariant object in $\Rep A$, letting $u^W_g: W \to W$,
be defined by $u^W_g(w) = g^{-1}w$, for every $g \in F$. We have
thus a functor $\mathcal G: \Rep (A\#_{\sigma}F) \to (\Rep A)^F$.
It is clear that $\mathcal F$ and $\mathcal G$ are inverse
equivalences of categories. This proves the proposition.
\end{proof}

We shall use later (c. f. Section \ref{clsfusionrules}) the
following description of irreducible representations of crossed
products given in \cite[Theorem 1.3]{MW}.

Consider  an irreducible $A$-module  $V$, and let $F^V$ its
stabilizer in $F$, and $\alpha$ the $2$-cocycle of $F^V$
associated to $V$ as in \cite[pp. 318]{MW}. Then there is an
equivalence between the category of $k_{\alpha}F^V$-modules and
the category of those $A\#_{\sigma}kF$-modules whose restriction
to $A$ is isomorphic to a direct sum of copies of conjugates of
$V$, which maps the $k_{\alpha}F^V$-module $W$ to the induced
module $\Ind_{A\#_{\sigma}kF^V}^{A\#_{\sigma}kF} U \otimes W$.

\subsection{Cocentral Hopf algebra extensions} Suppose that \begin{equation}\label{cleft}k \to A \to H
\overset{\pi}\to B \to k\end{equation} is a cleft exact sequence
of Hopf algebras; that is, the projection $\pi: H \to B$ admits a
convolution invertible $B$-colinear section $j: B \to H$. The
sequence \eqref{cleft} will be called \emph{cocentral} if
$\pi(h_1) \otimes h_2 = \pi(h_2) \otimes h_1$, for all $h \in H$.

If $B$ is finite dimensional, this is equivalent to saying that
the dual inclusion $\pi^*: B^* \to H^0$ is central. In this case
$B^*$ must be commutative, whence $B \simeq kF$ for some finite
group $F$.

We shall next give a description of the representation category of
a cocentral extension $H$ as in \eqref{cleft}, such that $B$ is
finite dimensional. That is, the Hopf algebra $H$ fits into a
cocentral cleft exact sequence
\begin{equation}\label{byagroup}k \to A \to H \to
kF \to k. \end{equation}  By the cleftness assumption, $H$ has the
structure of a bicrossed product $H \cong A {}^{\tau}\#_{\sigma}
kF$, with respect to a certain compatible datum $(., \rho, \sigma,
\tau)$, where $. : kF \otimes A \to A$ is a weak action, $\sigma:
kF \otimes kF \to A$ is an invertible cocycle, $\rho: kF \to kF
\otimes A$ is a weak coaction and $\tau: kF \to A \otimes A$ is an
invertible dual cocycle, subject to the compatibility conditions
in \cite[Theorem 2.20]{AD}.

\begin{lemma}\label{cocentral} The exact sequence \eqref{byagroup} is
cocentral if and only if the afforded weak coaction $\rho$ is
trivial. \end{lemma}

\begin{proof} The proof is straightforward. \end{proof}

In particular, $H$ is a crossed product $H \cong A \#_{\sigma} kF$
as an algebra. As in Subsection \ref{repcrossprod},  this gives
rise to an action of $F$ on $\Rep A$.

\begin{lemma} Suppose the exact sequence \eqref{byagroup} is
cocentral.  Then the afforded $F$-action on $\Rep A$ is by tensor
autoequivalences, and the functor $\mathcal F$ in Proposition
\ref{equiv} is a tensor functor.
\end{lemma}

\begin{proof} By Lemma \ref{cocentral} the afforded weak coaction is trivial.
On the other hand, by condition (2.21)  of compatibility between
the product and the coproduct in \cite{AD}, we have
$$\Delta_A(g.a) = \tau(g) (g.a_1 \otimes g.a_2)\tau(g)^{-1},$$ for
all $g \in F$, $a \in A$.  This implies that the action of
$\tau(g)^{-1}$ gives a well-defined
 $A$-linear isomorphism $\tau(g)^{-1}: T_g(V \otimes W) \to
T_g(V) \otimes T_g(W)$, natural in $V, W$, for all $g \in F$.
Moreover, the dual cocycle condition \cite[(2.18)]{AD} with
respect to the trivial coaction  $\rho$ implies that
$\tau(g)^{-1}$ actually gives a monoidal functor structure to the
autoequivalence $T_g$.

Let $g, h \in F$. By  condition (2.21) in \cite{AD} and the
assumption on $\rho$, we also get
$$\Delta(\sigma(g, h)) \tau(gh) = \tau(g)(g.\tau(h)) \; (\sigma(g,
h) \otimes \sigma(g, h)).$$ Inverting this identity, we find that
the map $f_{g, h}: T_g \circ T_h \to T_{gh}$ is indeed an
equivalence of tensor functors.

Formula \cite[(2.15)]{AD} for the comultiplication in $H$ reads in
this case \begin{equation}\label{coprod}\Delta(a \# g) =
\Delta(a)\tau(g)(g \otimes g),\end{equation} for all $a \in A$, $g
\in F$. This implies that $\mathcal F$ is a monoidal functor with
isomorphisms $\mathcal F(U \otimes V) \to \mathcal F(U) \otimes
\mathcal F(V)$ given by identities. This finishes the proof of the
lemma.
\end{proof}

In conclusion, we obtain the following description of the category
$\Rep H$. It generalizes the statement in \cite[Example 2.7]{nik}.

\begin{proposition}\label{equiv-hopf} Suppose that the cleft exact sequence \eqref{byagroup}
is cocentral. Then the functor $\mathcal F: (\Rep A)^F \to \Rep H$
defines an equivalence of tensor categories. \qed
\end{proposition}

\begin{remark} Suppose that $H$ is semisimple (hence finite dimensional).
Observe that the dual Hopf algebra $H^*$ fits into a
\emph{central} extension $k \to k^F \to H^* \to A^* \to k$. By
\cite[Proof of Theorem 3.8]{gel-nik}, this amounts exactly to the
fact that the category $\mathcal R = \Rep H^*$ is an $F$-graded
fusion category, $\mathcal R = \oplus_{g \in F} \mathcal R_g$,
with trivial component $\mathcal R_e = \Rep A^*$.

Therefore, Proposition \ref{equiv-hopf} establishes a duality
between $F$-graded fusion categories with trivial component $\Rep
A^*$ and $F$-equivariantizations of $\Rep A$, for a semisimple
Hopf algebra $A$. \end{remark}

\section{Extensions of group algebras by a group of order
2}\label{fusionrules}

Let $G$ be a finite group. In this section we shall assume that
$H$ is a semisimple Hopf algebra for which $G \subseteq G(H)$ as a
subgroup, and the group algebra $kG$ has index 2 in $H$. Therefore
$kG$ is a normal Hopf subalgebra of $H$ and we have an extension
\begin{equation}\label{extension}k \to kG \to H \to k\mathbb Z_2 \to k.
\end{equation}
This extension is necessarily cocentral, by \cite[Corollary
1.4.3]{Na1}. By Proposition \ref{equiv-hopf} there is an
equivalence of tensor categories
$$\Rep H \simeq (\Rep G)^{\mathbb Z_2},$$ with respect to an
appropriate action by tensor autoequivalences $\underline{\mathbb
Z_2} \to \underline{\Aut}_{\otimes}(\Rep G)$.

\begin{examples} (1) Nontrivial examples of Hopf algebras $H$ as in
\eqref{extension} are the duals of those whose representation
category is a Tambara-Yamagami category \cite{TY}, as observed in
\cite{gp-ttic}. In this case, the group $G$ is abelian.

\medbreak (2) Other examples of such $H$ appear in the
classification of Kac algebras of dimension $24$ by Izumi and
Kosaki \cite{IK}. More precisely, Kac algebras of dimension $24$
with group of group-likes of order $12$ are classified in
\cite[Theorem XIV.40-I]{IK}: for exactly two of the possible
isomorphism classes, the group of group-likes is nonabelian and
isomorphic to the alternating group $\mathbb A_4$.
\end{examples}

\medbreak Our first lemma describes the coalgebra structure of
$H$.

\begin{lemma}\label{coalgebra} Suppose $H$ is not cocommutative.
Then, as a coalgebra, $H \simeq kG \oplus \underset{n \text{
times}}{M_d(k)\oplus \dots \oplus M_d(k)}$, where $d > 1$ and $|G|
= d^2n$.
\end{lemma}

\begin{proof} As a coalgebra, $H \simeq kG \oplus M_{d_1}(k)^{n_1}\oplus \dots \oplus
M_{d_r}(k)^{n_r}$, where $d_i > 1$ and $n_i \geq 1$, $i = 1,
\dots, r$. Each subcoalgebra of the form $M_{d_1}(k)^{n_1}$ is
stable under left multiplication by $kG$ and is thus a left $(kG,
H)$-Hopf module. Therefore, by \cite{NZ}, $|G|$ divides
$n_id_i^2$, for all $i= 1, \dots, r$.

Then $\dim H = |G|+\sum_{i = 1}^rn_id_i^2 \geq (r+1)|G|$. Since,
on the other hand, $\dim H = 2|G|$, we get that $r = 1$, and $|G|
= d^2n$, $d = d_1$, $n = n_1$, as claimed. \end{proof}

\begin{corollary} Suppose  $|G|$ is square-free. Then $H$ is cocommutative.
 \qed \end{corollary}

The corollary implies that there is no semisimple Hopf algebra $H$
 with $[H: kG(H)] = 2$ and such that $|G(H)|$ is square-free.

\medbreak Let $X = \{\chi_1, \dots, \chi_n\}$ denote the set of
irreducible characters of degree $d$.

\medbreak Let $\chi$, $\psi$ be characters corresponding to the
$H$-comodules $V$ and $W$, respectively. Then the product
$\chi\psi$ is the character of the tensor product $V\otimes W$. We
may write $\chi \psi = \sum_{\mu} m(\mu, \chi\psi) \mu$,  where
$\mu$ runs over the set of irreducible characters and $m(\mu,
\chi\psi)$ denotes the multiplicity of $\mu$ in the product
$\chi\psi$.

Consider the action $G \times X \to X$ given by left
multiplication. Let $\chi \in X$ be an irreducible character, and
let  $G[\chi] \subseteq G$ be the stabilizer of $\chi$, that is,
$G[\chi]$ is the subgroup consisting of all  $g \in G$ such that
$g \chi = \chi$. In view of \cite[Theorem 10]{NR}, $g \in G[\chi]$
if and only if $m(g, \chi \chi^*) > 0$, if and only if $m(g, \chi
\chi^*) = 1$. Therefore
\begin{equation}\label{prod}\chi  \chi^* = \sum_{g \in G[\chi]} g +
\sum_{\epsilon(\mu) > 1} m(\mu, \chi \chi^*) \mu.\end{equation}
Note that $G[\chi^*] = \{ g \in G: \chi g = \chi \}$. On the other
hand, $G[\chi g] = G[\chi]$, $G[g\chi] = gG[\chi]g^{-1}$, for all
$g \in G$.

\begin{lemma}\label{transitiva} The group $G$ acts transitively on $X$ by left
multiplication. \end{lemma}

\begin{proof} Let $\chi \in X$ be the irreducible character corresponding to the simple subcoalgebra $C$.
We have $|G| = |G[\chi]| |G\chi|$, where $G\chi \subseteq X$
denotes the orbit of $\chi$. By \cite{NZ}, $|G[\chi]|$ divides
$d^2 = \dim C$, hence $|G| = d^2n$ divides $|G\chi| d^2$. This
implies that $n$ divides $|G\chi|$ and thus $|X| = n = |G\chi|$ as
claimed.
\end{proof}

\begin{corollary}\label{reglas} Let $\chi \in X$. Then $|G[\chi]| = d^2$.
In particular, $\chi\chi^* = \sum_{g\in G[\chi]}g$, for all $\chi
\in X$. \end{corollary}

\begin{proof} By Lemma \ref{transitiva}, $|G[\chi]| = d^2$ for all
$\chi \in X$. Combined with formula \eqref{prod} we get the
claimed expression for the product $\chi\chi^*$.
\end{proof}

\begin{remark}\label{right} Note that all arguments used so far can also be
applied to the right action of $G$ on $X$ given by right
multiplication. In particular, the right action is also
transitive.

As a consequence of this, we note the following. Let $\chi \in X$.
Then we may write $\chi^* = \chi g$, and therefore $G[\chi^*] =
G[\chi]$. This implies that $\chi\chi^* = \sum_{g\in G[\chi]}g =
\chi^*\chi$.
\end{remark}

\begin{corollary}\label{stabil} We have $G[\chi] = G[\psi]$, for all $\chi, \psi \in
X$. \end{corollary}

The common stabilizer $G[\chi]$, $\chi \in X$, will be denoted by
$\Gamma$.

\begin{proof} This follows from the fact that,  by Remark \ref{right}, the right action of
$G$ is also transitive. Thus $\psi = \chi g$ and $G[\psi] =
G[\chi]$. \end{proof}

\begin{proposition}\label{gamaab} The stabilizer $\Gamma$ is a normal
abelian subgroup of $G$ which admits a non-degenerate 2-cocycle.
\end{proposition}

Recall that a 2-cocycle $\alpha$ on $\Gamma$ is called
\emph{non-degenerate} if the twisted group algebra
$k_{\alpha}\Gamma$ is simple.

\begin{proof} Proposition 3.4.4 of \cite{Na1} implies that
$\Gamma$ is an abelian subgroup and admits a non-degenerate
2-cocycle.

Let $g \in G$ and fix $\chi \in X$. By Corollary \ref{stabil},
$G[g\chi] = G[\chi] = \Gamma$. On the other hand, $G[g\chi] =
gG[\chi]g^{-1} = g\Gamma g^{-1}$. Then $\Gamma$ is normal in $G$
as claimed.
\end{proof}

\begin{lemma}\label{ty} Let $\chi, \psi \in X$. Then $\psi \chi = \sum_{g \in \Gamma}
agb$, where $a, b \in G$ are such that $\psi = a\chi$ and $\chi =
\chi^*b$. In particular, $\psi \chi \in kG$, for all $\psi, \chi
\in X$.
\end{lemma}

\begin{proof} It is clear. \end{proof}

\begin{proposition}\label{abeliana} The group algebra $k\Gamma$ is a normal Hopf subalgebra of
$H$. There is an (abelian) exact sequence of Hopf algebras
\begin{equation*}k  \to k^{\widehat \Gamma} \to H \to kF \to k, \end{equation*}
where $F$ is a group of order $2n$. \end{proposition}

\begin{proof} Consider the group algebra $A = k\Gamma \simeq k^{\widehat \Gamma}$. It satisfies $AC = C =
CA$, for all simple subcoalgebras $C \subseteq H$, and also $\dim
A = \dim C$. By \cite[Proposition 3.2.6]{Na1} (see also Lemma
\ref{alternativo} in the Appendix), $A = k\Gamma$ is a normal Hopf
subalgebra of $H$.

By \cite[Corollary 3.3.2]{Na1}, the quotient Hopf algebra
$H/H(k\Gamma)^+$ is cocommutative, hence $H/H(k\Gamma)^+ \simeq
kF$ for some finite group $F$ of order $[H: k\Gamma] = 2n$.
\end{proof}

\begin{teor}\label{triangular} Let $A$ be a semisimple triangular Hopf algebra.
Let also $k  \to A \to K \to k\mathbb Z_2 \to k$ be an exact
sequence of Hopf algebras. Then $K$ is group-theoretical.
\end{teor}

\begin{proof} By the classification results of Etingof and Gelaki,
$A \simeq (kG)^J$ as a Hopf algebra, for some finite group $G$ and
$J \in kG\otimes kG$, a twist. Let $H = K^{J^{-1}}$, so that $H$
fits into an exact sequence $k \to kG \to H \to k\mathbb Z_2 \to
k$. By Proposition \ref{abeliana}, $H$ fits  into an abelian
extension as well. Hence $H$, and also $A$, must be
group-theoretical, in view of \cite{gp-ttic}. \end{proof}

\section{Classification}\label{clsfusionrules}

Let $H$ be a nontrivial Hopf algebra that fits into an exact
sequence \eqref{extension}. The aim of the present section is to
determine the structure of these Hopf algebras. We keep the
notation from previous sections.

\medbreak According to what we have proven in Section
\ref{fusionrules}, the isomorphism classes of simple objects in
the fusion category $\Rep H^*$ are parameterized by the set $G
\cup \{ x_1, \dots,  x_n\}$, where $n \geq 1$, $G$ is a group of
order $d^2n$, and $x_i$ are non invertible objects of dimension $d
> 1$, satisfying
\begin{equation}\label{tamyam}g \otimes h = gh, \quad x_i\otimes x_{i^*} = \oplus_{s \in \Gamma}s,
\quad x_i \otimes x_j \in kG, \end{equation} where $x_{i^*} =
x_i^*$, for all $g, h \in G$, $1 \leq i, j \leq n$.

\begin{lemma} Suppose, conversely, that $H$ is a noncocommutative
semisimple Hopf algebra such that $\Rep H^*$ has fusion rules as
described in \eqref{tamyam}. Then $G(H) \simeq G$ has index $2$ in
$H$ and $H$ fits into an exact sequence \eqref{extension}.
\end{lemma}

\begin{proof} The relation for the product $x_i\otimes x_{i^*}$ implies
that all $x_i$'s have the same dimension $d$. Hence, $\dim H =
|G|+d^2n$. Moreover, since $H$ is not cocommutative, then $d > 1$
and $G(H) \simeq G$.

Condition $x_i \otimes x_j \in kG$, for all $1\leq i, j\leq n$,
implies that the natural action of $G$ on the set $\{ x_1, \dots,
x_n\}$ is transitive; see \cite{NR}. Moreover, since $x_i\otimes
x_{i^*} = \oplus_{s \in \Gamma}s$, $\forall i$, then the
stabilizer of $x_i$ with respect to this action is the subgroup
$\Gamma$ and it has order $d^2$.

Therefore $|G| = d^2n$, implying that $G(H) \simeq G$ has index
$2$ in $H$ and thus $H$ fits into an exact sequence
\eqref{extension}.
\end{proof}

\medbreak Equations \eqref{tamyam} generalize the fusion rules
considered by Tambara and Yamagami in \cite{TY}, that correspond
in our setting to the case $n = 1$.

For the case $n=1$ the classification appears in the paper
\cite{tambara}, where, moreover, a necessary and sufficient
condition is given in order that a fusion category with the fusion
rules in \cite{TY}  admit a fiber functor, and hence be equivalent
to the representation category of some semisimple Hopf algebra.
The fact that, in this case, $H$ is an abelian extension as in
Proposition \ref{abeliana} was observed in \cite{gp-ttic}.

\medbreak By Proposition \ref{abeliana}, $H$ fits into an abelian
exact sequence
\begin{equation}\label{exact2}k  \to k^{\widehat \Gamma} \to H \to kF \to k, \end{equation}
where $\Gamma \simeq \widehat\Gamma$ is a normal abelian subgroup
of $G$ of order $d^2$ possessing a non-degenerate 2-cocycle, and
$F$ is a group of order $2n$.

In particular, as a Hopf algebra, $H$ is isomorphic to a bicrossed
product $H \simeq k^{\widehat\Gamma} {}^{\tau}\#_{\sigma}kF$,
corresponding to fixed actions $\tright: \widehat\Gamma \times F
\to F$, $\tleft: F \times \widehat\Gamma \to \widehat\Gamma$, and
compatible cocycles $\sigma: F \times F \to
(k^{\widehat\Gamma})^{\times}$, $\tau: \widehat\Gamma \times
\widehat\Gamma \to (k^{F})^{\times}$.

\medbreak We shall now compare the description of the
corepresentation theory of $H$ given by Lemma \ref{coalgebra} with
the corepresentation theory of crossed products given by Clifford
theory. See \cite[Section 3]{KMM}, \cite{MW}.

\medbreak Let $x\in F$ and let $\widehat\Gamma^x \subseteq
\widehat\Gamma$ denote the isotropy subgroup of $x$ with respect
to the action $\tright$.

For all $s, t \in \widehat\Gamma$, write $\tau(s, t) = \sum_{y\in
F} \tau_y(s, t)e_y$, where $e_y$, $y \in F$, are the canonical
idempotents in $k^F$, and  $\tau_y(s, t) \in k^{\times}$. The
restriction of $\tau_x$ defines a normalized 2-cocycle $\tau_x:
\widehat\Gamma^x \times \widehat\Gamma^x \to k^{\times}$. Let
$k_{\tau_x}\widehat\Gamma^x$ denote the corresponding twisted
group algebra.

\medbreak Since $H^* \simeq k^F \#_{\tau}k\widehat\Gamma$ is a
crossed product as an algebra,  the isomorphism classes of
irreducible $H^*$-modules (= $H$-comodules) are parameterized by
the modules
\begin{equation}\label{simple}V_{x, W} = \Ind_{k^F \#_{\tau}k\widehat\Gamma^x} \, x \otimes W
= H^* \otimes_{k^F \#_{\tau}k\widehat\Gamma^x} (x \otimes W),
\end{equation}
where $x$ runs over a set of representatives of the orbits of
$\widehat\Gamma$ on $F$,  and $W$ runs over a system of
representatives of isomorphism classes of irreducible left
$k_{\tau_x}\widehat\Gamma^x$-modules. We have $\dim V_{x, W} =
[\widehat\Gamma: \widehat\Gamma^x] \dim W$.

\begin{proposition}\label{actriv} The exact sequence \eqref{exact2} is cocentral.
In other words, the action $\tright: \widehat\Gamma \times F \to
F$ is trivial.  \end{proposition}

\begin{proof} Let $\mathcal C$ denote the fusion category of
representations of the dual Hopf algebra $H^*$, and view $\Rep
k^{\widehat\Gamma}$ as a fusion subcategory of $\mathcal C$. For
any simple object $x \in \mathcal C$ we have $xx^* = \sum_{s \in
\Gamma}s$. Hence $\mathcal C_{\ad} = \Rep k^{\widehat\Gamma}$,
where $\mathcal C_{\ad}$ is the adjoint subcategory as defined in
\cite[Section 8.5]{ENO}.

Let $U(\mathcal C)$ be the universal grading group of $\mathcal
C$, in the terminology of \cite[Section 3.2]{gel-nik}. Then the
category $\mathcal C$ is $U(\mathcal C)$-graded with trivial
component equal to $\mathcal C_{\ad}$.

By \cite[Theorem 3.8]{gel-nik}, $\mathcal C_{\ad} = \Rep
(H^*/H^*K^+)$, where $K = k^{U(\mathcal C)}$ is the maximal Hopf
subalgebra of $H^*$ which is contained in its center.  Exactness
of the sequence \eqref{exact2} now implies that $U(\mathcal C)$ is
isomorphic to $F$ and, moreover, the inclusion $k^F \subseteq H^*$
is central, as claimed.
\end{proof}

\begin{remark} Propositions \ref{equiv-hopf} and  \ref{actriv} imply that the category $\Rep H$ is
tensor equivalent to an $F$-equivariantization
$(\vect^{\Gamma})^F$. See Remark \ref{cocentralabeliana}.
\end{remark}

\medbreak Consider the subgroup $F_0 \subseteq F$ defined by $F_0
= \{ x \in F: \; [\tau_x] = 1 \}$.

\begin{proposition}\label{f0} There is an exact
sequence of groups $1 \to \Gamma \to G \to F_0 \to 1$. In
particular, the subgroup $F_0 \subseteq F$ is of order $n$.
\end{proposition}

\begin{proof} It follows from Lemma \ref{grouplikes}, since $F^{\widehat\Gamma} = F$.
By exactness of the sequence $1 \to \Gamma \to G \to F_0 \to 1$,
$|F_0| = n$ is a necessary and sufficient condition in order that
$[H: kG] = 2$.  \end{proof}

\begin{remark} Proposition \ref{f0} can be proved alternatively as follows.
Consider the simple $H^*$-modules \eqref{simple}. Note that $\dim
V_{x, W} = [\widehat\Gamma: \widehat\Gamma^x] \dim W = 1$, if and
only if $x \in F$ is a fixed point under the action of
$\widehat\Gamma$ and $W$ is an irreducible representation of
$k_{\tau_x}\widehat\Gamma = k_{\tau_x}\widehat\Gamma^x$ of
dimension $1$. If such $W$ exists, then the class of $\tau_x$ is
trivial (the twisted group algebra $k_{\tau_x}\widehat\Gamma$
being augmented), and therefore $k_{\tau_x}\widehat\Gamma \simeq
k\widehat\Gamma$.

Thus, since $\Gamma$ is abelian, every element $x\in F$ such that
$[\tau_x] = 1$ gives exactly $|\widehat\Gamma| = d^2$
one-dimensional $H$-comodules $V_{x, W}$'s. Hence $|G(H)| =
|F_0|d^2$. On the other hand, by Lemma \ref{coalgebra}, $H$ has
exactly $|G| = d^2n$ irreducible comodules of dimension 1.
Therefore $|F_0| = n$, as claimed. \end{remark}

\begin{proposition}\label{non-deg} For all $x \in F\backslash F_0$, the 2-cocycle $\tau_x \in H^2(\widehat\Gamma,
k^{\times})$ is nondegenerate. \end{proposition}

\begin{proof} Let $x \in F\backslash F_0$ and let $V_{x, W}$ be the simple $H$-comodule described in
\eqref{simple}, where $W$ is a simple $k_{\tau_x}\widehat\Gamma^x
= k_{\tau_x}\widehat\Gamma$-module. In view of Lemma
\ref{coalgebra} and Proposition \ref{f0}, we have $\dim V_{x, W} =
\dim W = d$. This implies that $\tau_x$ is nondegenerate, as
claimed. \end{proof}

We next prove our main classification result.

\begin{proof}[Proof of Theorem \ref{cls}] We need to show that a semisimple Hopf algebra
$H$ with $[H: kG(H)] = 2$ is determined by a triple $(\Gamma, F,
\xi)$, where $F$ is a finite group acting on the abelian group
$\Gamma$ by automorphisms, and $\xi = [(\sigma, \tau)] \in
\Opext(k^{\widehat \Gamma}, kF)$, such that the subgroup $F_0 : =
\{ x \in F:\, [\tau_x] = 1 \}$ has index $2$ in $F$, and $\tau_x$
is non-degenerate for all $x \in F\backslash F_0$.

\medbreak Suppose first given a semisimple Hopf algebra $H$ with
$[H: kG(H)] = 2$. Then $H$ fits into an exact sequence
\eqref{extension}, with $G = G(H)$. By the results in Section
\ref{fusionrules}, $H$ is an abelian extension $k \to k^{\widehat
\Gamma} \to H \to kF \to k$ associated to a matched pair
$(\widehat \Gamma, F)$ where the action $\widehat \Gamma \times F
\to F$ is trivial. Thus the action $\widehat \Gamma \times F \to
\widehat \Gamma$ is by group automorphisms, giving by
transposition an action by group automorphisms $F \to
\Aut(\Gamma)$.

Let  $\sigma: F \times F \to (k^{\widehat \Gamma})^{\times}$,
$\tau: \widehat \Gamma \times \widehat \Gamma \to
(k^{F})^{\times}$ be the associated $2$-cocycles. By Proposition
\ref{f0}, $[F: F_0] = 2$, and by Proposition \ref{non-deg},
$\tau_x$ is a non-degenerate $2$-cocycle, for all $x\in
F\backslash F_0$. In this way, we obtain a triple $(\Gamma, F,
\xi)$ satisfying the requirements. Furthermore, we have in this
case $G(H) \simeq \Gamma \rtimes_{\sigma} F_0$.

\medbreak Conversely, suppose given such a triple $(\Gamma, F,
\xi)$. Then we have a matched pair $(\widehat \Gamma, F)$ with
respect to the transpose action $\widehat \Gamma \times F \to
\widehat \Gamma$. These data give rise to a Hopf algebra $H =
k^{\widehat \Gamma} {}^{\tau}\#_{\sigma} kF$ as in \ref{ab-gral}.
By Lemma \ref{grouplikes} and the assumption $[F: F_0] = 2$, we
find that $[H: kG(H)] = 2$.

\medbreak Finally, suppose that the triples $(\Gamma, F, \xi)$,
$(\Gamma', F', \xi')$ give rise to isomorphic Hopf algebras $H$
and $H'$, respectively.

Let $(\sigma, \tau)$ be a pair of compatible cocycles representing
$\xi$. Since $[H: kG(H)] = 2$, then the coalgebra structure of $H$
is as described in Lemma \ref{coalgebra}. From Remark \ref{estab},
in view of the assumption that $\tau_x$ is non-degenerate for all
$x\notin F_0$, we get that $\Gamma \subseteq \Gamma_0$, where
$\Gamma_0$ is the (common) stabilizer in $G(H)$ of simple
subcoalgebras of dimension $>1$.

By construction, there is a cocentral exact sequence $k \to
k^{\widehat \Gamma} \to H \to kF \to k$. This implies that there
is a faithful $F$-grading on the category $\C$ of finite
dimensional $H$-comodules with trivial component $\C_e = k
\Gamma-\text{comod}$. Then necessarily $\C_{\ad} \subseteq
k\Gamma-\text{comod}$. See \cite{gel-nik}.

We have seen in the proof of Proposition \ref{actriv} that
$\C_{\ad} = k\Gamma_0-\text{comod}$. Hence $\Gamma_0 \subseteq
\Gamma$, and thus $\Gamma = \Gamma_0$. Similarly, $\Gamma'$
coincides with the stabilizer $\Gamma'_0$ of simple subcoalgebras
of dimension $> 1$ in $G(H')$.

Therefore, a Hopf algebra isomorphism $H \to H'$ must send
$\Gamma$ to $\Gamma'$. Then it induces an isomorphism of the
corresponding exact sequences. The rest of the theorem follows
from the fact that such exact sequences are classified by the
group $\Opext(k^{\widehat \Gamma}, kF)$.
\end{proof}

\begin{remark}\label{cls-tambara} Suppose that $F_0 = 1$. This corresponds to the
case $n = 1$ in Lemma \ref{coalgebra}; so that $\Rep H^*$ is a
Tambara-Yamagami category with $\Gamma$ as the group of invertible
objects.

Consider the invariant $(\Gamma, F, \xi)$ of $H$ given by Theorem
\ref{cls}, and identify $\Gamma \simeq \widehat \Gamma$. The group
$F$ is cyclic of order $2$, hence the action of $F$ on $\Gamma$
reduces to an automorphism $T\in \Aut (\Gamma)$ of order $2$.

Let $(\sigma, \tau)$ be a pair of cocycles representing $\xi$.
Since $F = \mathbb Z_2$, by \cite[Proposition 1.2.6]{pqq}, we may
assume that $\sigma = 1$. Thus $\xi$ reduces to a $2$-cocycle $\xi
= \tau_x: \Gamma \times \Gamma \to k^{\times}$, where $1\neq x \in
\mathbb Z_2$, such that $T^*(\xi) = \xi^{-1}$. Moreover, $\xi$ is
non-degenerate.

The non-degenerate symmetric bilinear form $\chi: \Gamma \times
\Gamma \to k^{\times}$ in \cite{TY} is given in this case by
$\chi(a, b) = \alpha(a, T(b))$, for all $a, b \in \Gamma$, where
$\alpha: \Gamma \times \Gamma \to k^{\times}$ is the
non-degenerate alternating form corresponding to $\xi$ under the
isomorphism $H^2(\Gamma, k^{\times}) \simeq \Hom(\Lambda^2\Gamma,
k^{\times})$.

Since the form defined by $\xi (T^*\xi_{21})^{-1}$ on $\Gamma$ is
symmetric, then there exists $\nu: \Gamma \to k^{\times}$ such
that $d\nu = \xi$.

We thus recover the invariants given by Tambara in
\cite[Proposition 3.2 and Theorem 3.5]{T}.
\end{remark}

\section{Semisimple Hopf algebras of low dimension}\label{lowdim}

In \cite{nik} a family of examples of non group-theoretical
semisimple Hopf algebras has been presented, answering a question
raised in \cite{ENO}. The smallest such example has dimension $36$
and it is a semisolvable Hopf algebra.

We shall show in this section that every semisimple Hopf algebra
of dimension $<36$ is group-theoretical, so that in fact $36$ is
the smallest possible dimension for that a non group-theoretical
semisimple Hopf algebra can have.

 As mentioned in the introduction, except for dimension
$24$, every semisimple Hopf algebra of dimension $<36$ is either
nilpotent (dimension $p$, $p^2$, $p^3$, $p^4$, $p^5$, where $p$ is
a prime number) or an abelian extension (dimensions $30$ and
$pq^2$, where $p$ and $q$ are prime numbers). Therefore all these
Hopf algebras are group-theoretical, in view of \cite{dgno,
gp-ttic}. We may then restrict our analysis to the case of
dimension $24$.

\medbreak In what follows we suppose that $H$ is a semisimple Hopf
algebra of dimension 24 over $k$.

\begin{proposition}\label{main}
Up to a cocycle twist of the multiplication or the
comultiplication, $H$ fits into an abelian extension $k \to
k^{\Gamma} \to H \to kF \to k$, where $|\Gamma||F| = 24$.
\end{proposition}

In particular, $H$ is group-theoretical, by \cite{gp-ttic}.

\begin{proof} We keep the notation and conventions in \cite[Lemma
6.2.1]{Na1}. By the results in \textit{loc. cit.}, $H$ fits into
an extension \begin{equation}\label{exacta}k \to A \to H \to
\overline H \to k,\end{equation} with $\dim A, \dim \overline H >
1$. We shall prove that, unless $H$ satisfies the claim,
$\overline H$ is necessarily cocommutative. Dualizing this
(observe that the statement in Proposition \ref{main} is of
self-dual nature), we also get that $A$ is necessarily
commutative, whence the exact sequence \eqref{exacta} can be
supposed to be abelian itself, and the claim is established.

Suppose $\overline H$ is not cocommutative. Then $\dim \overline H
= 6, 8$ or $12$.

\medbreak First of all, if $|G(H)| = 12$, that is, $[H: kG(H)] =
2$, then $H$ fits into an abelian extension, by Proposition
\ref{abeliana}. Hence we may assume that neither $H$ nor $H^*$
have group of group-likes of index 2.

Also, if $H$ contains a commutative normal Hopf subalgebra of
dimension 12, then $H$ is an abelian extension, since the quotient
has dimension 2.

Therefore we may assume that $H$ is not of type  $(1, 12; 2, 3)$
as a coalgebra. If $H$ is of type $(1, 3; 2, 3; 3, 1)$ as a
coalgebra, then it contains a commutative Hopf subalgebra of
dimension $12$, by \cite[Remark 2.2.2 (i)]{Na1}, and it is thus an
abelian extension; hence we may assume that $H$ is neither of this
type.

\medbreak If $\dim \overline H = 6$ and $\overline H$ is not
cocommutative, then $\overline H = k^{\mathbb S_3}$, and $6$
divides $|G(H^*)|$. Then we may assume that $H^*$ is of type $(1,
6; 3, 2)$ as coalgebra, by \cite[Lemma 6.1.1]{Na1}. In this case,
$H^*$ has a unique Hopf subalgebra $B$ of dimension 3, which
coincides with the group algebra of the stabilizer of a simple
subcoalgebra of dimension 9.

By \cite[Lemma 6.1.1 and Remark 6.1.2(ii)]{Na1}, $H$ has a Hopf
subalgebra $K$ of dimension 8 (since $\dim A = 4$ must divide
$|G(H)|$). Then necessarily $(H^*)^{\text{co} K^*} = B$ and $B$ is
normal in $H^*$. The coalgebra type of $H^*$ implies that
$H^*/H^*B^+$ is cocommutative \cite[Remark 3.2.7 and Corollary
3.3.2]{Na1}. Therefore, the extension $k \to B \to H^* \to K^* \to
k$ is an abelian extension.

\medbreak If $\dim \overline H = 8$, then $\dim A = 3$ and by
\cite[Lemma 6.1.1]{Na1} $H$ is of type $(1, 6;  3, 2)$ as a
coalgebra, since this is the only remaining possibility with
group-like elements of order 3.

In this case $G(H)$ has a unique (normal) subgroup $G$ of order
$3$, which must coincide with the stabilizer of all simple
subcoalgebras of dimension $9$, such that $A = kG$. As before, the
quotient $\overline H = H/H(kG)^+$ is cocommutative.

\medbreak Suppose that $\dim \overline H = 12$ and $\overline H$
is not cocommutative. If $\overline H$ is commutative, then $12$
divides $|G(H^*)|$ and we are done.

We may therefore assume that $\overline H$ is not cocommutative
and not commutative. By the classification in dimension 12
\cite{fukuda}, $\overline H$ is isomorphic to one of the self-dual
Hopf algebras $\mathcal A_0$ or $\mathcal A_1$, in the notation of
\cite[5.2]{Na1}.

If $\overline H \simeq \mathcal A_0$, then $\overline H$ is a
twisting of a group algebra; see \cite[Proposition 5.2.1]{Na1}.
Then $H^*$ is twist equivalent to a Hopf algebra with group of
group-likes of order 12, which must be an abelian extension.

\medbreak Suppose that $\overline H \simeq \mathcal A_1$. Then
$\overline H$ is of type $(1, 4; 2, 2)$ as a coalgebra and
$G(\overline H) \simeq \mathbb Z_4$ is cyclic of order $4$.  Also,
$\overline H \simeq \overline H^*$ is a Hopf subalgebra of $H^*$,
and thus $4$ divides $|G(H^*)|$.

Combining these with \cite[Lemmas 6.1.1 and 6.1.7]{Na1}, since we
may assume that $|G(H^*)| \neq 12$, we get that the possible
coalgebra types for $H^*$ are $(1, 4; 2, 5)$ and $(1, 8; 2, 4)$.

In the first case, $G(H^*) = G(\overline H^*)$ is cyclic. By
\cite[Proposition 2.1.3]{Na1}, $H^*$ contains a  Hopf subalgebra
$K$ of dimension $8$, which is not cocommutative and such that
$G(H^*) \subseteq G(K)$. But, since $K$ is not cocommutative, we
must have $G(K) \simeq \mathbb Z_2 \times \mathbb Z_2$, which is a
contradiction, because $G(H^*) \simeq \mathbb Z_4$. This discards
this possibility.

\medbreak Finally, consider the case where $H^*$ is of type $(1,
8; 2, 4)$ as a coalgebra. Again by dimension, we have $H^* =
k[\overline H^*, G(H^*)]$. Let $F \subseteq G(\overline H^*)$ be
the unique subgroup of order $2$. We have $F \subseteq Z(\overline
H^*)$. Also, $G(\overline H^*)$ is normal in $G(H^*)$. Therefore
$F$ is stable under the adjoint action of $G(H^*)$ and thus
central in $G(H^*)$. Then $F$ is central in $H^*$. Let $p: H^* \to
B$ be the projection to the quotient Hopf algebra $B =
H^*/H^*(kF)^+$. As before, $p(\overline H^*)$ and $p(kG(H^*))$ are
cocommutative and they generate $B$ as an algebra. Then $B$ is
cocommutative and $H$ is an abelian extension. This finishes the
proof of the proposition.
\end{proof}

\begin{remark} Suppose $k = \mathbb C$ is the field of complex numbers. It follows from Proposition \ref{main}
that, after twisting the multiplication or the comultiplication,
$H$ turns into a Kac algebra, therefore belonging to the list of
Izumi and Kosaki \cite[Chapter XIV]{IK}.\end{remark}

We have thus proved the following:

\begin{teor}\label{menor36} Let $H$ be a semisimple Hopf algebra of
dimension $< 36$. Then $H$ is group-theoretical.

In particular, $36$  is the smallest possible dimension that a non
group-theoretical semisimple Hopf algebra can have. \qed
\end{teor}

\section{Appendix: Characters and normality}\label{char}

Let $H$ be a semisimple Hopf algebra and let $A \subseteq H$ be a
Hopf subalgebra. In this Appendix we discuss some sufficient
conditions, in terms of the character multiplication of $H$, in
order that $A$ be normal in $H$.

\begin{lemma}\label{chardual} Suppose $\chi A\chi^* \subseteq A$, for all
irreducible characters $\chi \in H$. Then $A$ is normal in $H$.
\end{lemma}

\begin{proof}Let $C \subseteq H$ be the simple subcoalgebra with character
$\chi$; so that $\mathcal S(C)$ is the simple subcoalgebra with
character $\chi^*$.

Write $A = \bigoplus_{\lambda \in \Lambda}C_{\lambda}$, where
$\Lambda$ is the set of irreducible characters of $A$, and
$C_{\lambda} \subseteq A$ is a simple subcoalgebra with character
$\lambda \in \Lambda$. The character of $A$, as a left
$H$-comodule, is $\chi_A = \sum_{\lambda \in \Lambda}
\epsilon(\lambda)\lambda$.

By assumption, $\chi \lambda \chi^* \in A$ is a cocommutative
element, and therefore it can be written as a sum $\chi \lambda
\chi^* = \sum_{\lambda \in \Lambda} n_{\lambda}\lambda$, for some
$n_{\lambda} \geq 0$, $\lambda \in \Lambda$. This implies that,
for every $\lambda \in \Lambda$, $CC_{\lambda}\mathcal S(C)
\subseteq A$, because the multiplication map $H \otimes H \to H$
is a left $H$-comodule map.

Therefore, $CA\mathcal S(C) \subseteq A$, implying that for all
$c\in C$, $c_1A\mathcal S(c_2) \subseteq A$. Since the simple
subcoalgebras of $H$ span $H$, this proves that $A$ is normal.
\end{proof}

\begin{remark}\label{generadores} Suppose that  $C_1, \dots, C_r$ is a set of simple
subcoalgebras of $H$ that generate $H$ as an algebra. Then the
assumption $\chi A\chi^* \subseteq A$, for all irreducible
character $\chi \in H$, in Lemma \ref{chardual}, may be replaced
by $\chi_i A\chi_i^* \subseteq A$, for all $i = 1, \dots, r$,
where $\chi_i$ is the character of $C_i$. \end{remark}

\begin{lemma}\label{alternativo} Let $\chi \in H$ be an irreducible character corresponding
to the simple subcoalgebra $C \subseteq H$.

Suppose $G[\chi] = G[\chi^*]$ and $|G[\chi]| = \chi(1)^2$. Then
$kG[\chi]$ is normal in $k[C]$.
\end{lemma}

Here, $k[C] \subseteq H$ denotes the subalgebra generated by $C$,
which is a Hopf subalgebra of $H$.

\begin{proof} The assumption $|G[\chi]| = \chi(1)^2$
implies that $\chi\chi^* = \sum_{g\in G[\chi]}g$. On the other
hand, since $G[\chi] = G[\chi^*]$, we also have $\chi g = \chi$,
for all $g \in G[\chi]$. Then, for all $g \in G[\chi]$, $\chi g
\chi^* = \chi \chi^* \in kG[\chi]$. By Lemma \ref{chardual} and
Remark \ref{generadores}, $kG[\chi]$ is normal in $k[C]$.
\end{proof}

\bibliographystyle{amsalpha}

\begin{thebibliography}{A}

\bibitem{AD}   \textsc{N. Andruskiewitsch},  \textsc{J. Devoto},
\emph{Extensions of Hopf algebras},   St. Petersbg. Math. J.
\textbf{7} (1996),  17-52.

\bibitem{artamonov}   \textsc{V. A. Artamonov},
\emph{Semisimple finite-dimensional Hopf algebras}, Mat. Sb.
\textbf{198}  (2007), 3--28.

\bibitem{agaitsgory} {\sc S. Arkhipov} and {\sc D. Gaitsgory},
\emph{Another realization of the category of modules over the
small quantum group}, Adv. Math. \textbf{173} (2003), 114--143.

\bibitem{dgno} {\sc V. Drinfeld}, {\sc S. Gelaki}, {\sc D. Nikshych} and {\sc V. Ostrik},
\emph{Group-theoretical properties of nilpotent modular
categories}, preprint \texttt{arXiv:math.QA/0704.0195v2}.

\bibitem{eg-triang} {\sc P. Etingof} and {\sc S. Gelaki},
\emph{The classification of triangular semisimple and co-
semisimple Hopf algebras over an algebraically closed field}, Int.
Math. Res. Not. \textbf{2000} (2000),  223--234.

\bibitem{ENO}  {\sc P. Etingof}, {\sc D. Nikshych} and {\sc V. Ostrik},
\emph{On fusion categories}, Ann. Math. (2) \textbf{162}, 581-642
(2005).

\bibitem{fw}  {\sc E. Frenkel}, and {\sc E. Witten},
\emph{Geometric endoscopy and mirror symmetry}, Commun. Number
Theory Phys. \textbf{2}, 113--283 (2008).

\bibitem{fukuda} \textsc{N. Fukuda}, \emph{Semisimple Hopf algebras of dimension $12$},
Tsukuba J. Math. \textbf{21} (1997), 43--54.


\bibitem{gel-nik} {\sc S. Gelaki} and {\sc D. Nikshych}, \emph{Nilpotent fusion
categories}, Adv. Math. \textbf{217}  (2008), 1053--1071.

\bibitem{twist-simple} {\sc C. Galindo} and {\sc S. Natale}, \emph{Simple Hopf algebras and deformations of finite groups},
Math. Res. Lett. \textbf{14} (2007), 943-954.

\bibitem{IK}  \textsc{M. Izumi}  and  \textsc{H. Kosaki},
\emph{Kac algebras arising from composition of subfactors: general
theory and classification}, Mem. Amer. Math. Soc. \textbf{158},
750, (2002).

\bibitem{kassel}  {\sc C. Kassel},
\emph{Quantum groups},  Graduate Texts in Mathematics \textbf{155}
Springer-Verlag, New York (1995).

\bibitem{KMM}   \textsc{Y. Kashina},  \textsc{G. Mason}  and  \textsc{S. Montgomery},
\emph{Computing the Frobenius-Schur indicator for abelian
extensions of Hopf algebras},  J.  Algebra \textbf{251} (2002),
888--913.

\bibitem{ma-ext}  {\sc A. Masuoka},
\emph{Extensions of  Hopf algebras},
 Trabajos de Matem\' atica \textbf{41/99}, Universidad Nacional de C\' ordoba   (1999).

\bibitem{ma-ext2}  {\sc A. Masuoka},
\emph{Hopf algebra extensions and cohomology}, in: New Directions
in Hopf Algebras, MSRI Publ. \textbf{43} (2002), 167--209.

\bibitem{Mo} {\sc S. Montgomery},
 \emph{Hopf algebras and their actions on rings},  CMBS Reg. Conf. Ser. in Math. \textbf{82}, Amer. Math. Soc., 1993.

\bibitem{MW} {\sc S. Montgomery} and {\sc S. Whiterspoon},
 \emph{Irreducible representations of crossed products},  J. Pure Appl. Algebra {\bf 129} (1998), 315--326.

\bibitem{pqq}  \textsc{S. Natale},  \emph{On semisimple Hopf algebras of dimension $pq^2$},   J.  Algebra
\textbf{221} (1999), 242-278.

\bibitem{gp-ttic}  \textsc{S. Natale},  \emph{On group theoretical Hopf algebras
and exact factorizations of finite groups},   J.  Algebra
\textbf{270} (2003), 199-211.

\bibitem{Na1} {\sc S. Natale},
 \emph{Semisolvability of semisimple Hopf algebras of low dimension}, Mem. Amer. Math.
 Soc. \textbf{186} (2007).

\bibitem{NR} \textsc{W. Nichols}  and  \textsc{M. Richmond},
\emph{The Grothendieck group of a Hopf algebra}, J. Pure Appl.
Algebra \textbf{106} (1996), 297--306.

\bibitem{NZ}   \textsc{W. Nichols}  and  \textsc{M. Zoeller}, \emph{A Hopf
algebra freeness Theorem}, Amer. J.  Math. \textbf{111} (1989),
381--385.

\bibitem{nik}  \textsc{D. Nikshych},
\emph{Non group-theoretical semisimple Hopf algebras from group
actions on fusion categories},   preprint
\texttt{arXiv:0712.0585v1} (2007).

\bibitem{siehler1} \textsc{J. Siehler}, \emph{Braided near-group
categories}, preprint \texttt{arXiv:math/0011037v1}.

\bibitem{siehler2}  \textsc{J. Siehler}, \emph{Near-group categories}, Algebr. Geom.
Topol. \textbf{3} (2003), 719--775.

\bibitem{scha-bigalois}  \textsc{P. Schauenburg},
\emph{Hopf Bigalois extensions}, Commun. Algebra \textbf{24}
(1996) 3797--3825.

\bibitem{T}  \textsc{D. Tambara},
\emph{Representations of tensor categories with fusion rules of
self-duality for finite abelian groups}, Israel J. Math.
\textbf{118} (2000), 29--60.

\bibitem{tambara}  \textsc{D. Tambara},
\emph{Invariants and semi-direct products for finite group actions
on tensor categories}, J. Math. Soc. Japan \textbf{53} (2001),
429--456.

\bibitem{TY}  {\sc D. Tambara} and {\sc S. Yamagami},
\emph{Tensor categories with fusion rules of self-duality for
finite abelian groups}, J. Algebra \textbf{209} (1998), 692--707.

\end{thebibliography}

\end{document}